\documentclass[12pt]{article}
\usepackage{latexsym}
\begin{document}

\newcounter{thm}
\newtheorem{Def}[thm]{Definition}
\newtheorem{Thm}[thm]{Theorem}
\newtheorem{Lem}[thm]{Lemma}
\newtheorem{Cor}[thm]{Corollary}
\newtheorem{Prop}[thm]{Proposition}




                

\setlength{\baselineskip}{14pt}
\newcommand{\vlimsup}{\mathop{\overline{\lim}}}
\newcommand{\vliminf}{\mathop{\underline{\lim}}}
\newcommand{\Av}{\mathop{\mbox{Av}}}
\newcommand{\spec}{{\rm spec}}
\newcommand{\qed}{\hfill $\Box$ \\}
\newcommand{\subarray}[2]{\stackrel{\scriptstyle #1}{#2}}
\def\textmc{\rm}
\def\({(\!(}
\def\){)\!)}
\def\R{{\bf R}}
\def\Z{{\bf Z}}
\def\N{{\bf N}}
\def\C{{\bf C}}
\def\T{{\bf T}}
\def\E{{\bf E}}
\def\H{{\bf H}}
\def\Prob{{\bf P}}
\def\M{{\cal M}}     
\def\F{{\cal F}}
\def\G{{\cal G}}
\def\D{{\cal D}}
\def\X{{\cal X}}
\def\A{{\cal A}}
\def\B{{\cal B}}
\def\L{{\cal L}}
\def\a{\alpha}
\def\b{\beta}
\def\e{\varepsilon}
\def\de{\delta}
\def\ga{\gamma}
\def\k{\kappa}
\def\la{\lambda}
\def\fa{\varphi}
\def\th{\theta}
\def\si{\sigma}
\def\t{\tau}
\def\om{\omega}
\def\De{\Delta}
\def\Ga{\Gamma}
\def\La{\Lambda}
\def\Om{\Omega}
\def\Th{\Theta}
\def\lan{\langle}
\def\ran{\rangle}
\def\lbr{\left(}
\def\rbr{\right)}
\def\const{\;\operatorname{const}} 
\def\dist{\operatorname{dist}} 
\def\Tr{\operatorname{Tr}}
\def\quadd{\qquad\qquad}
\def\n{\noindent}
\def\beq{\begin{eqnarray*}}
\def\eeq{\end{eqnarray*}}
\def\supp{\mbox{supp}}
\def\beqn{\begin{equation}}
\def\eeqn{\end{equation}}
\def\bp{{\bf p}}
\def\sg{{\rm sgn\,}}
\def\1{{\bf 1}}

\begin{center}
{Asymptotic estimates of the distribution of Brownian  hitting time of a disc } \\
\vskip4mm
{K\^ohei UCHIYAMA} \\
\vskip2mm
{Department of Mathematics, Tokyo Institute of Technology} \\
{Oh-okayama, Meguro Tokyo 152-8551\\
e-mail: \,uchiyama@math.titech.ac.jp}
\end{center}

\vskip8mm
\n
{\it running head}:   Brownian  hitting time of a disc

\vskip2mm
\n
{\it key words}: Brownian hitting time,  asymptotic expansion, Fourier analysis, Laplace inversion
\vskip2mm
\n
{\it AMS Subject classification (2009)}: Primary 60J65,  Secondary 60G50.

\vskip6mm

\begin{abstract}
The distribution of the first hitting time  of a disc for the standard two dimensional Brownian motion is computed.    By investigating the inversion integral of its Laplace transform  we give  fairy detailed  asymptotic  estimates of its density valid uniformly with respect to the point where the Brownian motion starts from.  \end{abstract}
\vskip6mm

\section { Introduction}
\vskip2mm

Let $B_t$ be the standard two-dimensional Brownian motion and ${\bf t}_r^{(2)}$  the Brownian first hitting time  of  the disc of radius $r>0$ centered at the origin: ${\bf t}_{r}^{(2)}=\inf \{t>0: |B_t|\leq r\}$.  In this paper we obtain precise asymptotic estimates of the distribution function  
$P_x[{\bf t}_{r}^{(2)}>t]$
and its density
$$p_{r,x}(t)=P_x[{\bf t}_{r}^{(2)}\in dt]/dt,$$
where $P_x$ denotes the probability law of $(B_t)_{t\geq 0}$ started at $x$.
$p_{r,x}(t)$ is  invariant  under rotation of $x$ and the function $u(t,\xi)=p_{r, (\xi, 0)}(t)$ solves the diffusion equation $u_t=\frac12 u_{\xi\xi}+(2\xi)^{-1}u_\xi$ together with the initial-boundary conditions $u(0,\xi)=0$ $(\xi>r)$ and $u(\cdot, r)=\de$ (Dirac's delta function).
The   corresponding interior problem of estimating the exit  time distribution  from a circle, which admits a quite nice eigenfunction expansion,  is  well understood and completely different from the present problem.

An explicit form of the Laplace transform of $p_{r,x}$ is well known (see (\ref{T^(2)}) below) and its inversion gives
\beqn\label{q_r_t}
p_{r,x}(t)=\frac1{2\pi i}\int_{- i\infty}^{i \infty} \frac{K_0(|x|\sqrt{2z})}{K_0({r}\sqrt{2z})}e^{tz}dz,~~~~~~|x|>r.
\eeqn
Although  much  information for the distribution of ${\bf t}_{r}^{(2)}$ can be obtained directly from the property of the Laplace transform  itself  (an interesting  illustration of this  is found in \cite{IM}, problem 4.6.4, where the  limit  law of $\lg {\bf t}^{(2)}_r/\lg r^{-1}$ as  $r\to 0$ is found in a very simple way),  here we wish to find more or less explicit asymptotic expressions of  $p_{r,x}(t)$ as $t\to\infty$  (such an expression for $\int_0^tp_{r,x}(s)ds$ is derived  from  (\ref{q_r_t}) in \cite{S}, a weak version of  our Theorem \ref{thm3}, in the case when $x^2/t$ is constant). The motivation for the latter  comes from the corresponding problem for  the random walks for which the results of this paper are useful (see Remark 3 below).
\vskip2mm 

We derive two asymptotic estimates of $p_{r,x}(t)$, one  sharp when $|x|^2/ t$ is small and the other so when both $x$ and $t$ are large.
  For the former case the result is given in terms of  the function
\beqn\label{d_W}
W(\la)=\int_{0}^\infty  \frac{e^{- \la u}du}{[\lg u]^2+ \pi^2 }~~~~~~(\la>0).
\eeqn
We fix a positive constant ${r_\circ}$ once for all and put
$c_\circ=-2\ga + \lg (2/{r^2_\circ})$ so that $\frac12 r^2_\circ = e^{-(c_\circ+2\ga)}$, where
 $\ga=-\int_0^\infty (\lg u)e^{-u}du$ (Euler's constant). 
 
  We write
  $a\vee b= \max \{a,b\},  a\wedge b=\min\{a,b\}$, $\lg^+a=\lg (a \vee 1)$ for $a, b$ real and $x^2=|x|^2$.

  \vskip2mm      
 
\begin{Thm}\label{thm1} ~ Uniformly for $|x|>r_\circ$,  as $t\to\infty$
\beqn\label{BH_2}
p_{r_\circ, x}(t)=2\bigg [\lg \frac{|x|}{{r_\circ}}\bigg] e^{c_\circ}W(e^{c_\circ} t)+O\bigg(\frac{1+\lg ^+|x|}{t(\lg t)^2}\cdot\frac{x^2\wedge t}{t}\bigg).
\eeqn
\end{Thm}

\vskip2mm\n
 {\sc Remark} 1.~   The function $\la W(\la)$  admits the following asymptotic expansion in powers of $1/\lg \la$:
\beqn\label{Rmk2}
\la W(\la)=
\frac1{(\lg\la)^{2}}-\frac{2\ga}{(\lg \la)^{3}}- \frac{\textstyle{\frac12} \pi^2 -3\ga^2}{(\lg \la)^{4}}+\cdots
\eeqn
valid in the both limits as $\la\to\infty$ and as $\la \downarrow 0$ (cf. \cite{B}).  In Appendix we shall indicate how to derive this expansion; therein  one will also find three Fourier representations of $W$.  In \cite{U1} the present author uses the function
$N(\la):=\int_{\la}^\infty W(t)dt$
to express the leading term of the asymptotic formula of the mean number of visited sites of random walk on $\Z^2$; $N(\la)$  is called    the Ramanujan function in   \cite{W}  and   \cite{B}, where asymptotic expansions of a class of functions including it are obtained. 
 
\vskip2mm
For $|x|$  much smaller than $\sqrt t$, Theorem \ref{thm1} provides an effective asymptotic expansion of $p_{r_\circ}(t,x)$ in powers of $1/\lg (e^{c_\circ}t)$  in view of (\ref{Rmk2}),  while  in the region where $|x|$ is comparable with $\sqrt t$  the estimate of Theorem \ref{thm1} is poor.  In the latter case, i.e. the case when   $x^2/t$ is bounded away from both zero and  infinity,  the next theorem provides such an estimate that the ratio of  the error term to the leading term is not $O(1/(\lg t))$ but
 at most  $O(1 /(\lg t)^2)$, a  smaller order of magnitude than one that may be expected from Theorem \ref{thm1}. 
 \vskip2mm\n
 \begin{Thm}\label{thm2}~   Uniformly for $|x|>r_\circ$,  as $t\to\infty$
\begin{eqnarray}\label{q_x40}
p_{r_\circ, x}(t)=  \frac{\lg({\textstyle  \frac12} e^{c_\circ} x^2)\,}{t(\lg (e^{c_\circ} t))^2} e^{- x^2/2t} + \left\{\begin{array}{ll}   {\displaystyle \frac{2\ga \lg({t}/{x^2})}{t(\lg t)^3}+ O\bigg(\frac1{t(\lg t)^3} \bigg) }~~~~~ & \mbox{for}~~~x^2<t\\[4mm]
{\displaystyle  O\bigg(\frac{1+[\lg (x^2/t)]^2}{x^2(\lg t)^3}\, \bigg) }~~~~~& \mbox{for}~~~x^2 \geq t.
\end{array}\right.
\end{eqnarray}
\end{Thm}

\vskip2mm\n
 {\sc Remark} 2.~Substituting from the identity $1/r_\circ =(\frac12 e^{c_\circ})^{1/2}e^{\ga}$ into  the first term on the right side of  (\ref{BH_2})  one may realize that Theorem \ref{thm1} implies (\ref{q_x40}) if restricted to the region  $x^2\leq C t/(\lg t)^2$. 
The error estimate in (\ref{q_x40}) is best possible in any region where $x^2/t$ is bounded, while it can be improved significantly for  $x^2/t$ large (enough if larger roughly than  $4\lg\lg t$) (see  Remark 5 given at the end of Section 3).

\vskip2mm\n
 {\sc Remark} 3.~ For a random walk on the two-dimensional square lattice $\Z^2$ we have analogues of Theorems \ref{thm1} and \ref{thm2}  \cite{U2}.  Suppose that the walk is aperiodic and satisfies  $E X=0$, $E|X^{2+\de}]<\infty$ for some $\de\in [0,2].$  Let $a(x)$  denote  the potential function of the walk and $f_x(k)$  the probability that the random walk started from $x$ hits the origin for the first time at the $k$-th step. Then, uniformly in $x\in \Z^2$, as $k\to\infty$
$$f_x(k)=2\pi |Q|^{1/2}\,a^*(x)\,e^{c_\circ} {W(e^{c_\circ} k)}\Big[1+o\Big(k^{-\de/2}\Big)\Big]+ O\bigg(\frac{|x|_+^2}{k^2\lg k}\bigg),$$
where   $a^*(x)=\de_{0,x}+ a(x)$ and $c_\circ$ is a certain constant depending on the  walk.   As to Theorem {\ref{thm2} the very  analogue of (\ref{q_x40}) 
holds for $f_x(k)$ (one may simply replace $t$ by $k$ on its right side), provided  $\de>0$.

\vskip2mm
One can readily obtain an estimate of  $P_x[{\bf t}^{(2)}_{r_\circ}<t] =\int_0^t p_{r_\circ, x}(u)du$ from Theorem \ref{thm1} that is  relatively sharp  if $x^2< t/(\lg t)^2$ (see Remark 4 (iii) below), while it is  not so simple a matter  to find a proper asymptotic form of it from  Theorem \ref{thm2}. The latter one is given  in  the next theorem (sharp if $x^2> t/\lg t$), in which we use the following notation
$$\xi= \frac{|x|}{\sqrt t}, ~~~~~~  \fa (\a)= - \int_{1}^\infty \frac{e^{-\a y}}{y} \lg\bigg (1-\frac1{y}\bigg) \,dy~~~~(\a>0),$$ 
and
$$A_x(t)=  \frac{1}{\lg (e^{c_\circ} t)}\bigg[1-\frac{\ga}{\lg (e^{c_\circ} t)}\bigg]\int_{\xi^2/2}^\infty\frac{e^{-u}}{u}du +\frac{\fa(\xi^2/2)}{ [\lg (e^{c_\circ} t)]^2 }.$$ 

\vskip2mm
 \begin{Thm}\label{thm3}~ Uniformly for $|x|>r_\circ$,  as $t\to\infty$
\begin{eqnarray}\label{q_x42}
P_x[{\bf t}^{(2)}_{r_\circ}\leq t] 
= A_x(t)
+ \frac 1{(\lg t)^{3}} \times\left\{\begin{array}{ll} { O(\lg \frac12 \xi )}
~~~~~ & \mbox{for}~~x^2<t  \\[1mm]
{\displaystyle  O((\lg 2\xi)^2/\xi^2 )}~& \mbox{for}~~x^2 \geq t.
\end{array}\right.  
\end{eqnarray}
\end{Thm}

\vskip2mm\n
 {\sc Remark} 4.~(i)~ It holds that  $\fa(\a)=O(\a^{-1}e^{-\a} \log \a)$ as  $\a \to \infty$
 and $\fa(\a)= \frac16 \pi^2+ \a \lg \a +O(\a)$ as $\a\downarrow 0$.

 \vskip2mm
\n 
(ii)  On using the identity $\int_1^\infty e^{-u}u^{-1}du +\int_0^1(e^{-u}-1)u^{-1}du =-\ga$ 
$$\int_{\xi^2/2}^\infty\frac{e^{-u}}{u}du =\ga -\lg   (\xi^2/2) -\int_0^{\xi^2/2}\frac{e^{-u}-1}{u}du.$$
With the help of this  together with 
 $2\ga =\lg [2/e^{c_\circ} r_\circ^2] $ it is deduced that as $x^2/t\to 0$, 
$$A_x(t)=1- \frac{2\lg (|x|/r_\circ)}{ \lg (e^{c_\circ}t)  }  \bigg[1-\frac{\ga}{\lg (e^{c_\circ}t)}\bigg]  +  \frac{\frac16 \pi^2 - \ga^2}{ [\lg (e^{c_\circ}t)]^2 } (1+o(1) ) + \frac{\xi^2/2}{\lg (e^{c_\circ}t)}(1+O(\xi^2)),$$
which agrees with the expression of $P_x[{\bf t}^{(2)}_{r_\circ}\leq t] $  obtained from Theorem \ref{thm1} apart from the error  of magnitude $ O(|{\lg \frac12 \xi}|/{(\lg t)^{3}}) $ (see the next item).

 \vskip2mm
\n 
(iii) ~ Integrating the formula of Theorem \ref{thm1} leads to
$$P_x[{\bf t}^{(2)}_{r_\circ}>t] = \frac{2\lg (|x|/r_\circ)}{ \lg (e^{c_\circ}t)  }  \bigg[1-\frac{\ga}{\lg (e^{c_\circ}t)} - \frac{\frac16 \pi^2-\ga^2 }{ [\lg (e^{c_\circ}t)]^2 } + \cdots \bigg]  + O\bigg(\frac{\xi^2\lg |x|}{(\lg t)^2}\bigg) ~~~~~(x^2<t)  .$$
The error estimates  for $x^2<t$   in the formula  (\ref{q_x42}) and in this one     cannot be improved  since otherwise they had become inconsistent as is revealed by examining the sum of their principal terms in comparison to the  error terms. On equating   the  error terms   the latter formula  is sharper than the former  if  $\xi^2/\lg \xi=o( (\lg t)^{-2})$. 

 \vskip4mm
In the   dimensions $d\geq 2$ we have the formula
\beqn\label{T^(2)}
E_x[\exp\{-\la {\bf t}_r^{(d)}\}] =\frac{G_\la(|x|,r)}{G_\la(r,r) }=\frac{K_{d/2-1}(|x|\sqrt{2\la})|x|^{1-d/2}}{K_{d/2-1}(r\sqrt{2\la})r^{1-d/2}}~~~~~(\la>0)
\eeqn
 (cf. e.g., \cite{IM} \S 7.2, \cite{S}), where  $G_\la$ denotes the resolvent kernel for the $d$-dimensional Bessel process and  $K_\nu$  the  modified Bessel function of order $\nu$. 
If $d$ is odd, we have a more or less explicit expression of $p_{r, x}(t)$.  If $d$ is  even, any explicit expression of simple form can  not be expected,  but  an asymptotic  estimate of it with the error term whose ratio to $p_{r_\circ,x}(t)$ is  at most $O(\lg t /t)$ is easily obtained.  This is also true even if $d$ is not an integer (i.e., for the Bessel processes).  Finally  we note that
  Brownian scaling property gives the identity $p_{r/|x|, (1,0)}(t)= x^2 p_{r,x}(x^2 t)$,  or,  what is the same thing, 
$$P_{(1,0)}[{\bf t}_{r/|x|}^{(d)}>t] =
 P_x[{\bf t}_{r}^{(d)}> x^2t].$$
   
 We prove  Theorem \ref{thm1}   in Section 2 and those of Theorems \ref{thm2} and  \ref{thm3} in Section 3. In the last section some results on  $W(\la)$ are given.

\section{Proof of Theorem \ref{thm1}}
Put 
\beqn\label{g(u)}
g(z)=-\lg\, ({\textstyle \frac12} r_\circ \sqrt{2z}\,)-\ga;
\eeqn 
\beqn\label{K_0}
q_x(t)=\frac{1}{2\pi i}\int_{-i\infty}^{i\infty}\frac{K_0(|x|\sqrt{2z})}{g(z)}e^{tz}dz,
\eeqn
 and 
\beqn\label{K_00}
q^c_x(t)=\frac{1}{\pi}\int_{-\infty}^\infty\frac{K_0(|x|\sqrt{i2u})}{g(iu)}\cos tu\,du.
\eeqn
Here $i=\sqrt{-1}$; the integral in (\ref{K_0}) is along the imaginary axis;  $\lg z$ and $\sqrt z$ denote the principal branches in $-\pi <\arg x<\pi$  of the logarithm and the square root, respectively. 
\begin{Lem}\label{lem4.7} ~Uniformly for $|x|>r_\circ$, as $t\to\infty$
$$
p_{r_\circ, x}(t)-q_x(t) =O\bigg(\frac{1}{t^2(\lg t)}\wedge \frac1{|x|^4 (1+ \lg^+ |x|)}\bigg)
$$
and the difference $p_{r_\circ,x}(t)-q^c_x(t)$  admits the same estimate.
\end{Lem}
\vskip2mm\n
{\it Proof.}  We prove  the first relation only, the second one being proved in the same way. 
By definition 
\beqn\label{998}
K_0(z)=\sum_{k=0}^\infty \frac{(z/2)^{2k}}{(k!)^2}\bigg(\sum_{m=1}^k\frac1{m}-\ga-\lg(2^{-1}z)\bigg) ~~~~~~(-\pi<\arg z<\pi).
\eeqn
It follows that for $-1/2<u<1/2$
\beqn\label{g_u0}
K_0(r_\circ\sqrt{2iu}\,)=g(iu)+(r_\circ/2)^2 (2iu)(g(iu)+1)+O(u^2\lg |u|);
\eeqn
 hence
\beqn\label{g_u}
\frac1{K_0(r_\circ\sqrt{2iu}\,)}-\frac1{g(iu)}=(r_\circ/2)^2\frac{-2iu}{g(iu)}\bigg(1+\frac{1}{g(iu)}\bigg)+\frac{O(u^2)}{g(iu)}.
\eeqn
If  $h(z)$ is holomorphic  on $\Re z>0$,  then for any $\xi>0$
\beqn\label{G'}
\frac{d}{du} h(\xi\sqrt{i2u})= \frac1{2u}k(\xi\sqrt{i2u})~~~~\mbox{where}~~~k(z)=zh'(z).
\eeqn
In view of (\ref{g_u}), the identities  $K_0'(z)= - K_1(z)$ and $K'_{\nu+1}(z)=- \frac12[K_\nu(z)+K_{\nu+2}(z)]$ and the asymptotic formula 
\beqn\label{K0}
K_\nu(z)=(\pi/2z)^{1/2}e^{-z}(1+O(1/z))~~~~~~~~\mbox{as}~~~|z|\to\infty
\eeqn
$(-\pi<\arg z<\pi,~ \nu\geq 0)$ (cf. \cite{L} (5.11.9))  (as well as (\ref{q_r_t})), 
it therefore  suffices to show
\beqn\label{q(u,x)}
\int_{-\infty}^\infty\frac{2iu}{g(iu)}w(u)K_0(|x|\sqrt{2iu}\,)e^{itu}du=O\bigg(\frac{1}{t^2(\lg t)}\wedge \frac1{|x|^4 (1+ \lg^+ |x|) }\bigg),
\eeqn
where $w$ is a smooth function that equals 1 in a neighborhood of the origin  and vanishes outside a finite interval,  for  the  $1-w$ parts of the integrals in (\ref{q_r_t}) and (\ref{K_0})  are   $O(e^{-\e |x| } \wedge t^{-3})$ with $\e>0$, of which the  bound $O(t^{-3})$ is derived   by integrating by parts thrice.
In order to evaluate the integral  on the left side of (\ref{q(u,x)})  (as well as for the later use)  we bring in  the two functions
$$G(u)=G_x(u):= K_0(| x|\sqrt{2iu}\,)~~~~\mbox{and}~~~~F(z)= - z K_0'(z) =zK_1(z).$$
Noticing  that
$G'(u) =  -\, \frac1{2u} F(| x|\sqrt{2iu}\,)
$
and that   $F(z)$, $zF'(z)$ and $z^2 F''(z)$ are  all bounded  (see (\ref{FFF})  for   more precise estimates), we deduce that 
\beqn\label{F''}
|G(u)|\leq  C(1+ \lg^+ |u|^{-1}) ~ ~~  \mbox{and}~~~ | G^{(j)}(u)| \leq C/|u|^j~~   (j=1, 2, 3) 
\eeqn
valid uniformly for $|x|\geq r_\circ$.  
 Integrating by parts transforms the integral in (\ref{q(u,x)}) into
\beqn\label{F'''} -\, \frac{2}{t}\int_{-\infty}^\infty\frac{w(u)G(u)}{g(iu)}e^{itu}du - \frac{2}{t}\int_{-\infty}^\infty u\frac{d}{du}\bigg(\frac{w(u)G(u)}{g(iu)}\bigg)e^{itu}du.
\eeqn
In order to evaluate  the second integral  split its range of integral at $|u|=1/t$ and to the part $|u|>1/t$  apply integration by parts once or twice more, which with the help of   (\ref{F''}) gives  the bound $O(1/t^2\lg t)$  for the second term of  (\ref{F'''})  (cf. \cite{U2}, Lemma 2.2). For the first integral  in (\ref{F'''})  suppose $x^2< t$ and consider the function $\tilde G(u) := G(u) - g(iu)  - \log (r_\circ/|x|)$.   Then $\tilde G(u)$ is bounded  at least for  $|u|<1/\lg t$,  $\tilde G'(u)=O(1/u)$ and $\tilde G''(u)= O(1/u^2)$, and by the same procedure as  the one carried out right above  we obtain 
$$\int_{-\infty}^\infty \frac{w(u)\tilde G(u)}{g(iu) }e^{itu}du = O\bigg(\frac1{t\lg t}\bigg).$$
We must estimate the contribution of $g(iu)+\lg (r_\circ/|x|)$, but this is easily disposed of since the integral
$\int_{-\infty}^\infty [{-g(iu)}]^{-1}e^{itu}du$ equals  $AW(e^{c_\circ}t)$ for $t>0$ and 
$-Ae^{- e^{c_\circ}|t|}$ for  $t\leq 0$, where $A=4\pi e^{c_\circ}$  (see (\ref{W_1}) of Appendix), and hence is bounded by $C/|t| (\log t)^2$.
 
The other bound $O(1/|x|^4\lg( |x|+2) )$ is obtained by simply scaling the variable $u$ by $x^2$ and using   (\ref{K0}). The proof of the lemma is complete. \qed

  In view of Lemma \ref{lem4.7}   we have only to evaluate $q_x(t)$ for the proof of Theorem \ref{thm1}.  
  \begin{Lem}\label{lem4.8}~ Uniformly for $|x|> r_\circ $, as $t\to\infty$
  $$q_x(t)= 2e^{c_\circ}W(e^{c_\circ} t)\lg \frac{|x|}{r_\circ }+O\bigg(\frac{1+\lg^+|x|\,}{t( \lg t)^2}\bigg(\frac{x^2}{t}\wedge 1\bigg)\bigg).
$$  
  \end{Lem}
 \vskip2mm\n
 {\it Proof}.~ 
    Put $H(z)=K_0(\sqrt{2z})$, $-\pi<\arg z<\pi$, so that
$$q_x(t)=\frac1{2\pi i}\int_{-i\infty}^{i\infty} \frac{H(x^2z)}{g(z)}e^{tz}dz.$$ 
Since   $H(z)$ is analytic  and bounded by a constant multiple of $ -\lg ( |z| \wedge \frac12)$, Cauchy's  theorem gives
\beqn\label{q010}
q_x(t)= - \frac{1}{2\pi i}\int_{0}^\infty \frac{H(-x^2u +0i)}{g(-u+ 0i)}e^{-tu}du + \frac{1}{2\pi i}\int_{0}^\infty \frac{H(-x^2u-0i)}{g(-u-0i)}e^{-tu}du .
\eeqn
From the formula 
(\ref{998}) it follows that for $u>0$,
$$H(-x^2u\pm 0i)=\sum_{k=0}^\infty {(k!)^{-2}}{\Big(-\frac12 x^2u\Big)^{k}}\bigg(\sum_{m=1}^k\frac1{m} +g\Big( -(x/r_\circ)^2u\pm 0i\Big)\bigg),$$
which
combined with the identity  $g((x/r_\circ)^2z)=-\lg (|x|/r_\circ)+ g( z)$  gives
$$ \frac{H(-x^2u\pm 0i)}{g(-u\pm 0i)}= \frac{1}{g(-u\pm 0i)}\sum_{k=0}^\infty \frac{(-\frac12 x^2u)^{k}}{(k!)^2}\bigg(\sum_{m=1}^k\frac1{m}-\lg \frac{|x|}{r_\circ }\bigg) +\sum_{k=0}^\infty \frac{(-\frac12 x^2u)^{k}}{(k!)^2}.$$
Substitute this expression into the right side of (\ref{q010}). Noticing that  the contribution of the last infinite sum cancels out and that 
since $g(z)= -\frac12 \lg(e^{-c_\circ}z)$,
$$ \frac{-1}{i g(-u + 0i)}+  \frac{1}{i g(-u - 0i)} =\Im\frac{-2}{g(-u + 0i)} = \frac{-4\pi}{(\lg (e^{-c_\circ}u))^2+\pi^2},$$
we then deduce that
\beq
q_x(t)&=&\frac{1}{2\pi}\int_{0}^\infty  \frac{-4\pi}{(\lg (e^{-c_\circ}u))^2+\pi^2}\sum_{k=0}^\infty \frac{(-\frac12 x^2u)^{k}}{(k!)^2}\bigg(\sum_{m=1}^k\frac1{m}-\lg \frac{|x|}{r_\circ}\bigg) e^{-tu}du \\
&=& 2e^{c_\circ}W(e^{c_\circ} t)\lg \frac{|x|}{r_\circ }+\int_{0}^\infty \frac{-2e^{-tu}}{(\lg (e^{-c_\circ}u))^2+\pi^2}\sum_{k=1}^\infty \frac{(-\frac12 x^2u)^{k}}{(k!)^2}\bigg(\sum_{m=1}^k\frac1{m}-\lg \frac{|x|}{r_\circ }\bigg) du\\
&=& 2e^{c_\circ}W(e^{c_\circ} t)\lg \frac{|x|}{r_\circ }+O\bigg(\frac{(\lg( |x|+2))}{t( \lg t)^2}\bigg(\frac{x^2}{t}\wedge 1\bigg)\bigg)~~~~~~~~~~~~~~~(t\to\infty)
\eeq
uniformly for $|x|> r_\circ $ as desired. Here for the last equality we have applied the  crude bounds
 $$\sum_{k=1}^\infty\frac{(-y)^k}{(k!)^{2}}=O(y\wedge 1)~~~~\mbox{and}~~~~~\sum_{k=1}^\infty\frac{(-y)^k}{(k!)^{2}}\sum_{m=1}^k\frac1{m}=O(y\wedge 1)~~~~~~~~(y>0),$$
which are readily verified by using  $K_0(\sqrt{-4y})=O(e^{-2\sqrt y})$ ($y>1$) together with  the formulae (\ref{998})  and
  $$\sum_{k=0}^\infty\frac{(-y)^k}{(k!)^{2}}=J_0(2\sqrt y) =\frac{2}{\pi}\int_0^{\pi/2}\cos(2\sqrt y \sin \th)d\th  = O(y^{-1/4})~~ (y>1)$$
($J_0$ is the Bessel function of first kind of order 0). Thus  (\ref{BH_2}) has been verified.  \qed

\section{Proof of Theorems \ref{thm2} and \ref{thm3}}

 For the present purpose we  estimate the function  $q^c_x$ that is defined by (\ref{K_00}). 
 Theorem \ref{thm3} will be proved  by computing
$$Q_x(t): =\int_0^t q^c_x(s)ds =\frac1{\pi}\int_{-\infty}^\infty \frac{K_0(|x|\sqrt{i2u}\, )}{g(iu)}\cdot \frac{\sin tu}{u}du.$$
Owing to the fact noted in Remark 4 (iii)  we have only to consider the case $x^2> t/(\log t)^2$, when  the difference between $Q_x(t)$ and $\int_0^t p_{r_\circ, x}(u)du$ is   negligible in view of Lemma \ref{lem4.7}.   

The computation is based on the following  formula (\cite{L} (5.10.25)):
\beqn\label{*Q}
K_0(|x|\sqrt{i2u})=\frac12 \int_0^\infty \exp\bigg(-y- i\frac{x^2}{2y}u\bigg)\frac{dy}{y}.
\eeqn
 Let  $\xi =|x|/\sqrt t$ and 
$\fa(\a)=- \int_1^\infty {e^{-\a y}}{y^{-1}}\lg [1-y^{-1}]dy, ~ \a>0$ as in Introduction and define 
\beq
R(t) =  \frac{-2}{\pi [\lg (e^{c_\circ} t)]^2}\int_{-\infty}^\infty K_0(\xi \sqrt{i2u} \,)\frac{[\lg(iu)]^2}{\lg(iu/e^{c_\circ}  t)}\cdot \frac{\sin   u}{u}du.
\eeq
Note that $\fa(0+)= -\int_0^1 u^{-1}\lg(1-u) \, du=\pi^2/6$, in particular  $\fa$ is bounded.
\begin{Lem}\label{lem3.1} ~~ $Q_x(t)=A(t) + R(t)$, where
\beqn\label{5.0}
A(t) :=  \frac{1}{\lg (e^{c_\circ} t)}\bigg[1-\frac{\ga}{\lg (e^{c_\circ} t)}\bigg]\int_{\xi^2/2}^\infty\frac{e^{-y}}{y}dy + \frac{\fa(\xi^2/2)}{[\lg (e^{c_\circ} t)]^2}.
\eeqn
\end{Lem}

\vskip2mm\n
{\it Proof.}~ 
We let $c_\circ=0$ for simplicity and write $Q_x(t)$ in the form
$$Q_x(t) =\frac1{\pi}\int_{-\infty}^\infty \frac{K_0(\xi \sqrt{i2u}\, )}{g(iu/t)}\frac{\sin u}{u}du.$$
Since  $g(z)=- \frac12 \lg (z/e^{c_\circ})$, upon noting $\frac1{1-x}=1+x +\frac{x^2}{1-x}$
\beq
\frac1{2g(iu/t)}&=&\frac1{-\lg(iu/t)} =\frac1{\lg t}\cdot\frac1{1-{\lg(iu)}/{\lg t}}\\
&=&\frac1{\lg t}\bigg[1+\frac{\lg(iu)}{\lg t}\bigg] - \frac1{\lg (iu/t)}\bigg[\frac{\lg(iu)}{\lg t}\bigg]^2.
 \eeq
Hence 
\beqn\label{Qx}
Q_x(t) =\frac{I(t)}{\lg t}+\frac{II(t)}{(\lg t)^2}+ R(t),
\eeqn
where
$$I(t)=  \frac2{\pi }\int_{-\infty}^\infty K_0(\xi \sqrt{i2u} \,)  \frac{\sin u}{u}du  
 ~~~\mbox{and}~~~
II(t)=\frac2{\pi }\int_{-\infty}^\infty  \lg(iu)\, K_0(\xi \sqrt{i2u} \,)    \frac{\sin  u}{u}du. 
$$
Substitution from $(\ref{*Q})$ gives
$$I(t)= \frac{2}{\pi}\int_0^\infty \frac{\sin u}{u}du\int_0^\infty \frac{e^{-y}}{y}\cos \frac{\xi^2 u}{2y}dy.$$
Since
$$\int_0^\infty \frac{\sin u}{u}\, \cos \frac{\xi^2 u}{2y} du=\left\{\begin{array}{ll} \pi/2~~~~~~ &y>\frac12 \xi^2, \\[2mm]
0~~~~~~~&y<\frac12 \xi^2,
\end{array}\right.
$$
by interchanging the order of integration (easily justified) we obtain
$$I(t)=\int_{\xi^2/2}^\infty \frac{e^{-y}}{y}dy.
$$
Similarly 
$$II(t)= \frac{1}{\pi} \int_0^\infty \frac{e^{-y}}{y}dy \int_{-\infty}^\infty \frac{\sin u}{u}\exp\bigg(\frac{-i\xi^2}{2y} u\bigg) \lg (iu) du.$$
We  write  the inner integral on the right side  in the form
$$ \int_{-i \infty}^{i\infty} \frac{e^{(1-s)z}- e^{-(1+s)z}}{i2z} \lg z\, dz,~~~~~~~s=\frac{\xi^2}{2y}$$
and prove that this integral equals
 \beqn\label{C_int}
 \left\{\begin{array}{ll} 0~~~~~~~~&\mbox{if}~~s>1,\\
- \pi (\ga + \lg (1-s)) ~~~~~~~~&\mbox{if}~~s<1.
 \end{array}\right.
 \eeqn
To this end first consider the case $s<1$.  We apply  Cauchy's  theorem by considering   the term involving $e^{(1-s)z}$  in the left half plane and the other in the right half.   Note that  the integrand is discontinuous along the negative real line but regular otherwise.   Let $C^+(\e)$ and $C^-(\e)$ denote the right and left halves of the circle $z=\e e^{i\th}$ for small $\e>0$.   We then observe that the integral in question may be written as
  $$\bigg[-\int_{-\infty+0i}^{-\e+0i}+\int_{-\infty-0i}^{-\e -0i}+\int_{C^-(\e)}\bigg] \frac{e^{(1-s)z}}{i2z} \lg z \,dz + \int_{C^+(\e)} \frac{e^{-(1+s)z}}{i2z} \lg z \,dz$$
 Noting that  $\lg (-u\pm 0i)=\lg u \pm i\pi$ for $u>0$   we rewrite this in the form
$$\pi \int_\e^\infty \frac{e^{-(1-s)u}}{u}\log u\,du + \frac12 \int_{-\pi}^\pi (\lg \e+i\th) d\th +O(\e \lg \e).$$
On letting $\e\downarrow 0$ with the help of  the identity $\lim_{\e\downarrow 0}[\int_\e^{\infty}e^{-au}u^{-1}du+\lg \e] = - \ga -\lg a$ ($a>0$) this is reduced to the second expression of  (\ref{C_int}). The case $s>1$ would  now be  obvious.

Now
one can conclude that
\begin{eqnarray}\label{II}
II(t) &=& -\ga \int_{\xi^2/2}^\infty \frac{e^{-y}}{y}dy - \int_{\xi^2/2}^\infty \frac{e^{-y}}{y}\lg\frac{y-\frac12 \xi^2}{y}dy \nonumber \\
&=& -\ga \int_{\xi^2/2}^\infty \frac{e^{-y}}{y}dy + \fa( \xi^2/2),
\end{eqnarray}
and substitution into (\ref{Qx}) completes the proof of the lemma.\qed

\n
{\it Proof of Theorem  \ref{thm3}}.~
We must evaluate $R(t)$.   We claim that
   \beqn\label{R21} 
R(t) =\left\{\begin{array}{ll} O(|\lg \frac12 \xi| /{(\lg t)^3})~~~~ &\mbox{if} ~~~~\xi<1,\\
O((\lg (2\xi))^2/\xi^2(\lg t)^3) ~~~~ &\mbox{if} ~~~~\xi\geq 1,
\end{array}\right.
\eeqn
which is enough for  the proof of   Theorem \ref{thm3}.

First consider the case $\xi  <1$. We split the integral defining $R(t)$ at $|u|= 1$ and $|u| =1/\xi^2$. For simplicity we consider only the integral on $u>0$ and let $c_\circ=0$.  Let $D(1)=(0,1), D(2)=(1,1/\xi^2], D(3)= (1/\xi^2,\infty)$ and put
\beqn\label{lemB}
B_j= \int_{D(j)} K_0(\xi\sqrt{i2u} \,)\frac{[\, \lg(iu)]^2}{\lg (iu/t)}\cdot \frac{\sin u}{u}du ~~~~(j=1,2,3). 
\eeqn
In view of the  bound $K_0(z) =O(\lg \frac14 |z|)$  ($|z|<1$) we  immediately obtain  that  
$|B_1| =O(|\lg \frac12 \xi| /{\lg t)}.$
For $B_2$   perform integration by parts  (with $\sin u$ to be  integrated).  Using the identity (\ref{G'}) with $|x|$ replaced by $\xi$  we then  infer that
\begin{eqnarray}\label{B2}
|B_2| &\leq& \frac{C}{\lg t} +C\int_1^{1/\xi^2}\bigg[  \frac{\,|F(\xi\sqrt{i2u}\,)| \cdot|\lg (iu)|^2}{|\lg (iu/t)|} \\
&&~~~~~~~~~~~~~~~~~~ + |K_0(\xi\sqrt{i2u})| \frac{\,  ( |\lg (iu/t)|+1) \cdot |\lg (iu)|^2 \,}{|\lg (iu/t)|^2}
\bigg] \frac{du}{u^2}.\nonumber
 \end{eqnarray}
Since $F(z)=- zK_1(z)=O(1)$  ($|z|<1$), this shows that  $|B_2|  =O(|\lg \frac12 \xi| /{\lg t}).$  On changing the variable of integration  
 \beqn\label{B3}
 B_3=  \int_1^\infty K_0( \sqrt{i2u} \,)\frac{[\, \lg(iu/\xi^2)]^2}{\lg (iu/|x|)}\cdot \frac{\sin u/\xi^2}{u}du.
 \eeqn
 The integration by parts then gives $|B_3|\leq C \xi^2|\lg \xi|^2/ |\lg |x|  =O(|\lg \frac12 \xi| /{\lg t}).$ Thus the claim is verified in the case $\xi<1$.
 
For the case $\xi\geq 1$  we make a similar argument but without integration by parts. Employing  the bound given  in  (\ref{K0})  for $|u|\geq 1/\xi^2$, we obtain
$$|R(t)|\leq \frac{C}{|\lg t|^2}\Bigg[ \int_0^{1/\xi^2} \frac{|\lg(\frac12 \xi\sqrt u) | \cdot |\lg(iu)|^2}{|\lg (iu / t)|}du +\int_{1/\xi^2}^\infty \frac{e^{- \xi\sqrt{u}}|\lg(iu)|^2}{(\xi \sqrt u\,)^{1/2}|\lg (iu/t)|}\cdot \frac{| \sin u|}{|u|}du   \Bigg]$$ 
The first integral in the square brackets and the part $\int_{1}^\infty$ of the second one are easily evaluated to be  $O((\lg \xi)^2/\xi^2\lg t)$.
 For the part $\int_{1/\xi^2}^1$  we  first   replace $|\lg (iu/t)|$ in the denominator by $\lg t$ and then make change of the variable by $y=\xi^2 u$ to find it  to be   $O((\lg \xi)^2/\xi^2\lg t)$. Hence $R(t)= O((\lg \xi)^2/\xi^2(\lg t)^3)$, as required.
\qed

We are to derive Theorem \ref{thm2}  from Theorem \ref{thm3} by  evaluating the derivatives $A'(t)$ and  $R'(t)$.
\begin{Lem}\label{lem3.2} ~~ There exists a constant $C$ such that for $t>2$
$|R'(t)|\leq  C/t (\lg t)^3$  for $x^2<t$ and  $|R'(t)|\leq  C[1+(\lg \xi)^2]/|x|^2 (\lg t)^3$   for $x^2 \geq t$.
\end{Lem}
 \vskip2mm\n
{\it Proof}.~  As before let $c_\circ=0$ for simplicity.
Differentiate the defining expression of $R(t)$  and observe that in the integrand there then appears 
$$K_0'(\xi \sqrt{i2u} \,)\sqrt{i2u}\, \xi' =-K_0'(\xi \sqrt{i2u} \,)\sqrt{i2u}\,\xi /2t .$$ 
 With the function  $F(z)=- zK_0'(z)=zK_1(z)$, we can write  the result as
\begin{eqnarray}\label{R'}
R'(t)&=&  \frac{-1}{\pi [\,\lg  t]^2t}\int_{-\infty}^\infty F(\xi \sqrt{i2u} \,)\frac{[\,\lg(iu)]^2}{\lg(iu/ t)}\,\frac{\sin u}{u}du  \nonumber \\
&&+  \frac{-2}{\pi [\,\lg  t]^2t}\int_{-\infty}^\infty K_0(\xi \sqrt{i2u} \,)\frac{[\,\lg(iu)]^2}{[\,\lg(iu/ t)]^2}\,\frac{\sin u}{u}du  \nonumber\\
&&+  \frac{4}{\pi [\,\lg  t]^3t}\int_{-\infty}^\infty K_0(\xi \sqrt{i2u} \,)\frac{[\,\lg(iu)]^2}{\lg(iu/ t)}\,\frac{\sin u}{u}du.
\end{eqnarray}

Consider the case $\xi<1$. By the same argument that derives the bound of $B_j$ in the proof of Theorem  \ref{thm3}  the last two terms are both $O((\lg \xi) /t(\lg t)^4)$.  It holds that 
\beqn\label{FFF}
F(z) =1+ O(z^2\lg |z|)~~~\mbox{ and}~~ 
F'(z)=z\lg z + O(z)~~~~~( |z|<1) 
\eeqn
 (as being readily deduced from (\ref{998}));   and 
$F'(z) = K_1(z) -2^{-1} z[K_0(z)+ K_2(z)] =O(|z|^{1/2} e^{-z})$      ($ |z| >1/2$). 
Obviously the first  integral in (\ref{R'})  restricted on $|u|<1$ is $O(1/\log t)$.  The integral on $1<u<1/\xi^2$  is also $O(1/\log t)$ as is deduced in the same way as $B_2$ is estimated  (one has only to replace  in (\ref{B2})   $K_0$ and $F$  by $F$ and $zF'$, respectively).  Similarly the integral on $1/\xi^2<u<\infty$ is evaluated to be at most $O(1/\log t)$ (one may replace $K_0$ by $F$ in (\ref{B3})).
Hence  $R'(t)=O(1/(\lg t)^3t)$ as $t  \to\infty$  uniformly for $\xi <1$.

In the case $\xi\geq 1$ one uses the bound $F(\xi\sqrt{i2 u})=O(e^{- \xi \sqrt {|u|/2}})$ to see that the integral involving $F$ is at most a constant multiple of 
$$\frac1{\lg t}\int_{0}^1 e^{-\frac12 \xi s}(\lg s)^2 sds =O\bigg(\frac{1+(\lg \xi)^2}{\xi ^2\lg t}\bigg),$$
the integral on $|u|>1$ being much smaller than this  for  large $\xi$. 
The other integrals are estimated in a similar way.  \qed

\n
{\it Proof of Theorem 2}.  Recalling  $\ga=- \int_0^\infty e^{-y}\lg y dy$, observe that 
$$\fa'(\a)=- (\ga +\lg \a)\frac{e^{-\a}}{\a} - \frac1{\a}\int_\a^\infty \frac{e^{- y}}{y}dy.$$
Then  some easy computation leads to
\begin{eqnarray} \label{B'1}
A'(t) = \frac{\lg( {\textstyle \frac12} e^{c_\circ} x^2)}{t[\lg (e^{c_\circ} t)]^2}e^{- x^2/2t}+
  \frac{2\ga}{t [\lg (e^{c_\circ} t)]^3}\int_{\xi^2/2}^\infty\frac{e^{-u}}{u}du - \frac{2\fa(\xi^2/2)}{t [\lg (e^{c_\circ} t)]^3}.
\end{eqnarray}
where  $A(t)$ is defined by (\ref{5.0}). Considering the case $x^2<t$ 
 let $b(t)$ denote the sum of the first two terms on the right side of (\ref{q_x40}), namely
$$b(t)= \frac{\lg( {\textstyle \frac12} e^{c_\circ} x^2)}{t[\lg (e^{c_\circ} t)]^2}e^{- x^2/2t} + \frac{2\ga \lg^+({t}/{x^2})}{t(\lg t)^3},$$
 to obtain
\begin{eqnarray} \label{B'}
A'(t) = b(t)+O\bigg(\frac{1}{t(\lg t)^3}  \bigg).
\end{eqnarray}
  Combined with Lemma \ref{lem3.2} this shows Theorem \ref{thm2} in the case $\xi <1$. The case $\xi\geq 1$  is  also  deduced  from  (\ref{B'1}) 
 in a similar way.  \qed

\n
{\sc Remark 5.} ~(i)~ Let $c_\circ=0$.  Then computations similar  to those performed for derivation of  (\ref{II}) lead to 
\beqn\label{last}
 q_x(t) 
=  \int_1^\infty  \frac{e^{-\frac12 \xi^2/ y}}{y} e^{-(y-1)t}dy
  - \int_{0}^1 \frac{e^{- \xi^2/2y}}{y}  dy\int_0^\infty \frac{e^{-(1-y) t u}}{(\lg u)^2 +\pi^2}du
 \eeqn
and
$$Q_x(t)  =    \int_{0}^1 \frac{e^{-\xi^2/2y}}{y}dy\int_0^\infty \frac{e^{-(1-y)u}du}{([\lg (u/t)]^2 +\pi^2)u} +S_x(t) $$
with $S_x(t)=  - \int_1^{\infty}\frac{e^{-\xi^2/2y}}{y}e^{-(y -1)t}dy
+ e^{-t} \int_0^\infty\frac{e^{-y}}{y}e^{-x^2/2y}dy= O(1/t)$ uniformly in $|x| > r_\circ$.
(In the derivation of (\ref{last}) we need to interchange the order of the repeated integral that arises as one substitutes  (\ref{*Q}) into the defining expression of $q_x(t)$. The argument for justification of  the interchange  involves somewhat delicate analysis.)
Although it  seems hard  to derive   Theorems \ref{thm2} or \ref{thm3}  from these expressions,  they are useful if $\xi$ is large: at least we can derive from (\ref{last}) that uniformly for $|x| \lg |x| < t< x^2$, as $|x| \to\infty$ 
$$q_x(t) =\frac{e^{-x^2/2t}}{t \lg t}(1+o(1)).$$
(ii)~  If one uses the representation 
$$K_0(| x|\sqrt{2iu}\,)=\frac1{\pi}\int_0^\pi d\a\int_0^\infty \frac{\cos (|x| r \sin \a)}{2iu+r^2}rdr
$$
 (\cite{E}, p.45 (15))) in place of (\ref{*Q}),  certain computations using  Cauchy's  theorem  lead to a result similar to Theorem \ref{thm2}, but this way is  more  involved than that we have adopted.

\section{Appendix}

(A)~  In order to derive the expansion (\ref{Rmk2})  observe
$$\la W(\la) =\frac{-1}{\pi} \Im\int_0^\infty \frac{e^{-\la u}\la du}{\lg u + i\pi}= \frac1{\pi \lg \la} \int_0^\infty \Im\bigg(1- \frac{\lg u+i\pi}{\lg \la}\bigg)^{-1}e^{-u}du.$$
On  using   the formula $\int_0^\infty (\lg u)^2 e^{-u}du= \frac16 \pi^2+\ga^2$ (\cite{E} p.169 (13))  a standard argument then leads to   (\ref{Rmk2}).  By  a more sophisticated method Bouwkamp \cite{B} derives an asymptotic expansion
for a class of functions defined by similar Laplace transforms; as a special case  it gives that
$\la W(\la)= \sum_{n=0}^\infty {c_n}{ (\log \la)^{-n-1 }},$ with $(c_n)$ determined by 
$$ \sum_{n=0}^\infty \frac{c_n}{n!}z^n  = \frac{z}{\Ga(1-z)}.$$
(Although in \cite{B} only the case $\la \to\infty$  is considered, the method is valid for the case $\la \downarrow 0$.)

\vskip4mm\n
(B)~ Here we derive   three Fourier representations of $W(\la)$.  The first one is
\beqn\label{W_1}
\frac1{2\pi}\int_{-\infty}^\infty\frac{e^{-i\la u}}{\lg(-iu)}\,du=\left\{\begin{array}{ll}W(\la)~~~~&(\la>0), \\[1mm]
-e^\la~~~~&(\la<0).
\end{array}\right.
\eeqn
For $\la>0$  this is obtained  by  Cauchy's  theorem with the help of the equality $\lg(-(u\pm 0i))$ $=\lg u\mp i\pi$ valid for $u>0$.
On the other hand,   by calculus of residues we see that for  $\la<0$, 
$$ \frac1{2\pi}\int_{-\infty}^\infty  \frac{e^{-i \la t}\,dt}{\lg(-it)} =\frac{-1}{2\pi}\int_{|z-1|=1}\frac{e^{\la z}}{\lg z}\frac{dz}{i}=- e^{\la}.$$
It is noted that since $\lg (-it)=\lg|t|-i\frac12 \pi\,\sg t$, the last identity means  
$$
\frac12\int_{0}^\infty  \frac{\sin |\la| t}{[\,\lg t]^2+\frac14 \pi^2 }\,dt=\frac1{\pi}\int_{0}^\infty  \frac{\,\lg t}{[\,\lg t]^2+\frac14 \pi^2 }\,\cos \la t\,dt+e^{-|\la|},
$$
which together with (\ref{W_1}) leads to 
\beqn\label{W_2}
W(\la)=\int_{0}^\infty  \frac{\sin \la t}{[\,\lg t]^2+\frac14 \pi^2 }\,dt -e^{-\la} =\frac2{\pi}\int_{0}^\infty  \frac{(\lg t)\, \cos \la t}{[\,\lg t]^2+\frac14 \pi^2 }\,dt+e^{-\la}~~~~~~(\la>0).
\eeqn

\vskip4mm
{\bf Acknowledgments.}~ The author wishes to thank the anonymous  referee for his carefully reading  the original manuscript,  providing several valuable comments, and  in particular pointing out the paper \cite{W} and, through it, \cite{B}.

\end{document}